\documentclass[sn-mathphys]{sn-jnl}

\jyear{2022}%

\theoremstyle{thmstyleone}%
\theoremstyle{thmstyletwo}%

\theoremstyle{thmstylethree}%

\begin{document}

\title[Solving the Cauchy problem for~a~three-dimensional difference equation]{Solving the Cauchy problem for~a~three-dimensional difference equation in a parallelepiped}


\author[1,2]{\fnm{Marina} \sur{Apanovich}}\email{marina.apanovich@list.ru}
\equalcont{These authors contributed equally to this work.}

\author*[2,3]{\fnm{Alexander} \sur{Lyapin}}\email{aplyapin@sfu-kras.ru}
\equalcont{These authors contributed equally to this work.}

\author[1]{\fnm{Konstantin} \sur{Shadrin}}\email{kvsh\_buffon@mail.ru}
\equalcont{These authors contributed equally to this work.}

\affil*[1]{\orgname{Kransnoyarsk State Medical University}, \orgaddress{\city{Krasnoyarsk}, \postcode{660022}, \country{Russia}}}

\affil[2]{\orgname{Siberian Federal University}, \orgaddress{\city{Krasnoyarsk}, \postcode{660041}, \country{Russia}}}

\affil[3]{\orgname{Fairmont State University}, \orgaddress{\city{Fairmont}, \postcode{26554}, \state{WV}, \country{United States}}}


\abstract{The aim of this article is further development of the theory of linear difference equations with constant coefficients.  We present a new algorithm for calculating the solution to the Cauchy problem for a three-dimensional difference equation with constant coefficients in a parallelepiped at the point using the coefficients of the difference equation and Cauchy data. The implemented algorithm is the next significant achievement in a series of articles justifying the Apanovich and Leinartas' theorems about the solvability and well-posedness of the Cauchy problem. We also use methods of computer algebra since the three-dimensional case usually demands extended calculations.}




\maketitle

\section{Introduction}\label{sec1}

Difference equations arise in various fields of mathematics and have numerous applications in science and technology. For example, in combination with the method of generating functions they appear to be a powerful method of study enumerative problems in the combinatorial analysis (see, for example, \cite{ref-book1, ref-book2, ref-journal1}). Another source of difference equations is discretization of differential equations. For example, the discretization of the Cauchy-Riemann equation resulted in the theory of discrete analytical functions (see \cite{ref-journal2, ref-journal3}), which is used in the theory of Riemann surfaces and combinatorial analysis (see, for example, \cite{ref-journal4, ref-journal5}). Methods of discretization of a differential problem are an important part of the theory of difference schemes and also lead to difference equations (see, for example, \cite{ref-book3, ref-book4}). Also difference equations are used in discrete time dynamical models  (see \cite{ref-book5, ref-book6}).

Additional conditions ("initial", "boundary", "Cauchy data") are given for the space of solutions of a multidimensional difference equation, which allows to select a unique solution from an infinite set of solutions, and the corresponding problem is called the Cauchy problem for a difference equation. A one-dimensional case is simple since the inital data of the Cauchy problem is finite (thus a gererating function for the solution is rational). The multidimensional case is different since the initial data is given on an infinite set and might have a difficult structure. A significant number of works is devoted to the development of algorithms for solving multidimensional difference equations with constant and polynomial coefficients (see, for example, \cite{ref-journal6, ref-journal7, ref-journal8}), and in \cite{ref-journal8a} difference equations with coefficients in a form of rational functions are considered. The connection between the generating function for solution to the Cauchy problem for a two-dimensional difference equation with constant coefficients and the generating function of the initial data is studied in \cite{ref-journal9, ref-journal10}. In \cite{ref-journal11} an algorithm for calculating the solution of the Cauchy problem for a two-dimensional difference equation with constant coefficients at a point from the coefficients of the difference equation and the initial data of the Cauchy problem was developed and implemented.

However, the multidimensional cases have not been investigated properly. The three-dimensional case is of great importance in problems of thermodynamics (thermal conductivity, see \cite{ref-book3}, anisotropic diffusion, see \cite{murrayjd}), mathematical biology (distribution of morphogen in insect wings \cite{murrayjd}, distribution of population density species \cite{bratus}). This paper considers the first step -- solving the Cauchy problem for a three-dimensional difference equation with constant coefficients in a parallelepiped at a point from coefficients of the difference equation and the initial data of the Cauchy problem. We give a notation of the Cauchy problem in a "parallelepiped" and describe a function of the initial data on this set in a way that allows us to develop and implement a computer algebra algorithm.

\section{The Cauchy problem for a polynomial difference operator in a parallelipiped}

Let $f$ be a function of integer variables $f:\mathbb{Z}_{+}^{n}\rightarrow\mathbb{C}$ and $\delta_{j}$ be a shift operator by $j$-th variable, i.e. $\delta_{j}f\left(t_{1},...,t_{n}\right)=f\left(t_{1},...,t_{j-1},t_{j}+1,t_{j+1},...,t_{n}\right)$, $j=1,2,...,n$.
If $\alpha=\left(\alpha_{1},...,\alpha_{n}\right)$ is a multi-index, then $\|\alpha\|=\alpha_{1}+...+\alpha_{n}$,
$\delta^{\alpha}=\delta_{1}^{\alpha_{1}}...\delta_{n}^{\alpha_{n}}$.
For two multi-indexes $\alpha=\left(\alpha_{1},...,\alpha_{n}\right)$ and $\beta=\left(\beta_{1},...,\beta_{n}\right)$
the inequality $\alpha\leqslant\beta$
means that $\alpha_{j}\leqslant\beta_{j}$ for all $j=1,2,...,n$.

We consider a difference polynomial operator
\begin{equation}\label{Lyapin:01}
P\left(\delta\right)=\sum_{0\leqslant\alpha\leqslant m}c_{\alpha}\delta^{\alpha},
\end{equation}
where
$c_{\alpha}$ are constants coefficients of the operator $P\left(\delta\right)$.
We denote $\alpha=\left(\alpha_{1},...,\alpha_{n}\right)$ and $'\alpha=\left(\alpha_{1},...,\alpha_{n-1}\right)$,
and then write the operator $P\left(\delta\right)$
as follows
\begin{equation}
P\left(\delta\right)=\sum_{\alpha_{n}=0}^{m_{n}}\left(\sum_{'0\leqslant'\alpha\leqslant'm}c_{'\alpha,\alpha_{n}}\delta^{'\alpha}\right)\delta_{n}^{\alpha_{n}}.
\end{equation}

We call a polynomial $$P\left('z,z_{n}\right)=\sum_{\alpha_{n}=0}^{m_{n}} \left(\sum_{'0\leqslant'\alpha\leqslant'm}^{}c_{'\alpha,\alpha_{n}}'z^{'\alpha}\right) z_{n}^{\alpha_{n}}=\sum_{\alpha_{n}=0}^{m_{n}}P_{j}\left('z\right)z_{n}^{\alpha_{n}}$$
as a characteristic polynomial for the difference operator $P\left(\delta\right)$, and its degree by a variable $z_{n}$ is an order of the difference operator $P\left(\delta\right)$.

Let $x=\left(x_{1},...,x_{n-1},0\right)$ be a point in the integer lattice $\mathbb{Z}_{+}^{n}$
and $$\Pi_{'x}=\left\{ t\in\mathbb{R}_{+}^{n}:'0\leqslant't\leqslant'x\right\}$$ be parallelepiped of a dimension $n-1$ in the hyperplane $t_{n}=0$.

We fix a point
$\beta=\left(\beta_{1},...,\beta_{n-1},m_{n}\right)\in\mathbb{Z}_{+}^{n}$, $'0\leqslant'\beta\leqslant'm$
such that $c_{\beta}\neq0$,
and consider a set
$\Pi_{'\beta,'m}=\left\{ 't\in\Pi_{'x}:'\beta\leqslant't\leqslant'x-'m+'\beta\right\} $.
Let
$L=\left(\Pi_{'x}\setminus\Pi_{'\beta,'m}\right)\times\left[0,x_{n}\right]$
be a set where the initial data is given. We consider a problem: \textit{find a solution to the difference equation}

\begin{equation}\label{Lyapin:02}
P\left(\delta\right)f\left(x\right)=g\left(x\right),x\in\Pi=\Pi_{'x}\times\left[0,x_{n}\right],
\end{equation}
\textit{such that the condition}
\begin{equation}\label{Lyapin:03}
f\left(x\right)=\varphi\left(x\right),x\in L,
\end{equation}
is valid, where
$g\left(x\right)$
and
$\varphi\left(x\right)$
are given functions of integer variables.

The problem (\ref{Lyapin:02}) -- (\ref{Lyapin:03}) is called the Cauchy problem for a polynomial difference operator (\ref{Lyapin:01}).

In \cite{ref-book7} the stability of a homogeneous two-layer linear difference scheme with constant coefficients is investigated for $m_{n}=1$. In \cite{ref-journal12} the solvability of the problem \ref{Lyapin:02} -- \ref{Lyapin:03} was investigated for $n=2$. In the theory of difference schemes, such problems are multilayer implicit difference scheme. In \cite{ref-journal13} the well-posedness of problem (\ref{Lyapin:02}) -- (\ref{Lyapin:03}) for $n = 2$ is investigated and an easily verifiable sufficient condition for correctness is proved. In \cite{ref-journal14} for $n = 3$, an easily verified sufficient condition for the solvability of the Cauchy problem (\ref{Lyapin:02}) -- (\ref{Lyapin:03}) is proved.

We consider Cauchy problem (\ref{Lyapin:02}) -- (\ref{Lyapin:03}) for $n=3$ and denote $\mathbb{Z}^{3}$ --
the integer lattice  and $\mathbb{Z}^{3}_{+}$ is a subset of this lattice consisting of points with non-negative integer coordinates. Let $\delta_{1}, \delta_{2}, \delta_{3}$ be shift operators by variable $x, y, z$ accordingly, which means that $\delta_{1}f(x,y,z)=f(x+1,y,z)$, $\delta_{2}f(x,y,z)=f(x,y+1,z)$, $\delta_{3}f(x,y,z)=f(x,y,z+1)$. We denote a parallelepiped $\Pi=\left\{ (x,y,z)\in\mathbb{Z}_{+}^{3},\,0\leqslant x\leqslant B_{x},\,0\leqslant y\leqslant B_{y},\,z \geqslant 0\right\}$ in a positive octant of the integer lattice $\mathbb{Z}^{3}_{+}$, and $B_{x}+1$ is a width of the parallelipiped $\Pi$ and $B_{y}+1$ is a length of the parallelipiped $\Pi$. The difference polynomial operatop (\ref{Lyapin:01}) is

\begin{equation}\label{Lyapin:04}
P(\delta_{1},\delta_{2},\delta_{3})=\sum_{j=0}^{m}\sum_{i_{y}=0}^{b_{y}}\sum_{i_{x}=0}^{b_{x}}c_{i_{x}i_{y}j}
\delta_{1}^{i_{x}}\delta_{2}^{i_{y}}\delta_{3}^{j}=\sum_{j=0}^{m}P_{j}(\delta_{1},\delta_{2})\delta_{3}^{j},
\end{equation}
where $b_x$, $b_y$ and $m$ define a size of a difference scheme and
\begin{equation}
P_{j}(\delta_{1},\delta_{2})=\sum_{i_{x}=0}^{b_{x}} \sum_{i_{y}=0}^{b_{y}}c_{i_{x}i_{y}} \delta_{1}^{i_{x}}\delta_{2}^{i_{y}}, j=0,1,...,m.
\end{equation}

The characteristic polynomial is
\begin{equation}
P(s,w,v)=\sum_{j=0}^{m}\sum_{i_{x}=0}^{b_{x}}\sum_{i_{y}=0}^{b_{y}}c_{i_{x}i_{y}j}s^{i_{x}}w^{i_{y}}v^{j},
\end{equation}
where $m$ is an order of the difference operator $P(\delta_{1},\delta_{2},\delta_{3})$, $b_{x}<B_{x}$, $b_{y}<B_{y}$.

We fix $\beta=\left(x_{\beta},y_{\beta},m\right)$ such that $c_{x_{\beta}y_{\beta}m}\neq0$
and consider a set $\Pi_{\beta}=\{(x,y,z)\in\mathbb{Z}_{+}^{3}:\,0\leqslant x-x_{\beta}\leqslant B_{x}-b_{x},\,0\leqslant y-y_{\beta}\leqslant B_{y}-b_{y},\,z>m-1\}$,
then $L=\Pi\setminus\Pi_{\beta}$ and the Cauchy problem is

\textit{to find a solution to the difference equation}
\begin{equation}\label{Lyapin:05}
P(\delta_{1},\delta_{2},\delta_{3})f(x,y,z)=g(x,y,z),(x,y,z)\in\ensuremath{\Pi}
\end{equation}
\textit{under a condition that}
\begin{equation}\label{Lyapin:06}
f(x,y,z)=\varphi(x,y,z),(x,y,z)\in L,
\end{equation}

where $g(x,y,z)$ and $\varphi(x,y,z)$ are given functions of integer arguments.

Problem (\ref{Lyapin:05}) -- (\ref{Lyapin:06}) is the Cauchy problem for polynomial difference operator (\ref{Lyapin:04}).

In \cite{ref-journal14} it was proven that problem (\ref{Lyapin:05}) -- (\ref{Lyapin:06}) is uniquely solvable if the condition
\begin{equation}\label{Lyapin:06a}
\|c_{x_{\beta}y_{\beta}m}\| > \sum_{\left(\alpha_{1}, \alpha_{2}\right) \neq \left(x_{\beta}, y_{\beta}\right), \alpha_{3} = m} \|c_{\alpha_{1}\alpha_{2}\alpha_{3}}\|
\end{equation}
takes place.

Our problem is to compute a value of function $f\left(x,y,z\right)$ at point $A$ with coordinates $\left(x_{1},y_{1},z_{1}\right)$.

\section{Description of the Input Data}

A solution to Cauchy problem (\ref{Lyapin:05}) -- (\ref{Lyapin:06}) for three-dimensional difference equation with constant coefficients at a point $A$ with coordinates $(x_{1},y_{1},z_{1})$ is a value of the function $f(x,y,z)$ at a point $A$. The algorithm of computing a value of function $f(x,y,z)$ at a point with given coordinates $(x_{1},y_{1},z_{1})$ is recursive and reduces to computing values of a function $f(x,y,z)$ on the finite set of points $\left(x,y,z\right)$ in the set $L=\Pi\setminus\Pi_{\beta}$.

Initial data (6) is giver by a matrix $F$ of a dimension three, containing a finite set of values of the initial data of the Cauchy problem. Coefficients of three-dimensional difference equation are given by a matrix $C$ of dimension three. For the technical implementation of the algorithm, matrices of coefficients $C$ and initial data $F$ are specified in layers, starting from the lowest one. Let us illustrate the procedure for specifying the matrices $F$ and $C$.

For a difference equation
\begin{equation}\label{Lyapin:07}
\begin{array}{c}
c_{000}f\left(x,y,z\right)+c_{100}f\left(x+1,y,z\right)+c_{200}f\left(x+2,y,z\right)+c_{010}f\left(x,y+1,z\right)+\\
+c_{110}f\left(x+1,y+1,z\right)+c_{210}f\left(x+2,y+1,z\right)++c_{001}f\left(x,y,z+1\right)+\\
+c_{101}f\left(x+1,y,z+1\right)+c_{201}f\left(x+2,y,z+1\right)+c_{011}f\left(x,y+1,z+1\right)+\\
+c_{111}f\left(x+1,y+1,z+1\right)+c_{211}f\left(x+2,y+1,z+1\right)=0
\end{array}
\end{equation}

the first layer of the matrix $C$ of coefficients is a matrix $C_{0}$:

\begin{equation}
C_{0} =
\begin{pmatrix}
c_{000} & c_{100} & c_{200}\\
c_{010} & c_{110} & c_{210}
\end{pmatrix},
\end{equation}

the second layer of $C$ is a matrix $C_{1}$:

\begin{equation}
C_{1} =
\begin{pmatrix}
c_{001} & c_{101} & c_{201}\\
c_{011} & c_{111} & c_{211}
\end{pmatrix}.
\end{equation}

A matrix $C$ of coefficients is written as follows:
\begin{equation}
C =
\left(\begin{array}{c|c}
\begin{array}{c} C_0 \end{array} & \begin{array}{c} C_1 \end{array}
\end{array}\right) =
\left(
\begin{array}{c|c}
\begin{matrix}
c_{001} & c_{101} & c_{201}\\
c_{011} & c_{111} & c_{211}
\end{matrix}
&
\begin{matrix}
c_{001} & c_{101} & c_{201}\\
c_{011} & c_{111} & c_{211}
\end{matrix}
\end{array}
\right)
\end{equation}

For difference equation (\ref{Lyapin:07}) we have $b_{x}=2$, $b_{y}=1$, $m=1$. For $\beta=\left(1,1,1\right)$, $B_{x}=4$, $B_{y}=2$, the three-dimensional matrix $F$ of the initial data is
$F =
\left(\begin{array}{c|c}
\begin{array}{c} F_0 \end{array} & \begin{array}{c} F_1 \end{array}
\end{array}\right)
$,
where
\begin{equation}
F_{0}=\left(\begin{array}{ccccc}
\varphi\left(0,0,0\right) & \varphi\left(1,0,0\right) & \varphi\left(2,0,0\right) & \varphi\left(3,0,0\right) & \varphi\left(4,0,0\right)\\
\varphi\left(0,1,0\right) & \varphi\left(1,1,0\right) & \varphi\left(2,1,0\right) & \varphi\left(3,1,0\right) & \varphi\left(4,1,0\right)\\
\varphi\left(0,2,0\right) & \varphi\left(1,2,0\right) & \varphi\left(2,2,0\right) & \varphi\left(3,2,0\right) & \varphi\left(4,2,0\right)
\end{array}\right),
\end{equation}
\begin{equation}
F_{1}=\left(\begin{array}{ccccc}
\varphi\left(0,0,1\right) & \varphi\left(1,0,1\right) & \varphi\left(2,0,1\right) & \varphi\left(3,0,1\right) & \varphi\left(4,0,1\right)\\
\varphi\left(0,1,1\right) & * & * & * & \varphi\left(4,1,1\right)\\
\varphi\left(0,2,1\right) & * & * & * & \varphi\left(4,2,1\right)
\end{array}\right).
\end{equation}

The entries of $F$, denoted by $*$, are calculated when the algorithm is executed. However, it is not possible to calculate the element $\varphi (1,2,1)$ without calculating the elements $\varphi \left (2,2,1 \right)$, $ \varphi (1,1,1)$, $\varphi (2,1,1)$, $\varphi (3,1,1)$, $\varphi (3,2,1)$. Thus, to find the unknown elements, it is necessary to solve a system of linear difference equations of the form (8) using the initial data $\ensuremath{\varphi(i,j,k)}$,
where $i=0,...,4$, $j=0,1,2$, $k=0$, $\varphi(i,j,k)$, $i=0,...,4$, $j=0$, $k=1$ and $\varphi\left(0,1,1\right)$, $\varphi\left(0,2,1\right)$, $\varphi\left(4,1,1\right)$, $\varphi\left(4,2,1\right)$.

Finally, the input data is finite:
\begin{enumerate}
  \item a three-dimensional $(b_{x}+1)\times(b_{y}+1)\times(m+1)$-matrix $\ensuremath{C=(c_{\alpha_{1}\alpha_{2}\alpha_{3}})}, \ensuremath{\alpha_{1}=0, \ldots, b_{x}}, \ensuremath{\alpha_{2}=0, \ldots,b_{y}}$, $\alpha_{3}=0,...,m$ with coefficients $\ensuremath{c_{\alpha_{1} \alpha_{2} \alpha_{3}}}$ of three-dimensional difference equation;
  \item a point $\beta_{0}=(x_{\beta},y_{\beta},m)$;
  \item a point $A$ with coordinates $(x_{1}, y_{1}, z_{1})$, which defines the coordinates of the desired value of the function $f(x, y, z)$ and a number of layers in the three-dimensional matrix $F$ of initial data;
  \item a three-dimensional $(B_{x}+1)\times(B_{y}+1)\times(z_1+1)$-matrix of the initial data $F=(\varphi(x,y,z))$ for $(x,y,z)\in L$ and $\varphi(x,y,z)=0$ for $(x,y,z)\notin L$.
\end{enumerate}

Since the coordinates of the elements of the three-dimensional matrix $C$ of coefficients of the difference operator and the three-dimensional matrix $F$ of initial data in the Cartesian coordinate system $(X, Y, Z)$ do not coincide with their coordinates in the matrix (row $ \times $ column $ \times $ layer), then they have to be transformed from Cartesian coordinates
$\left(D\left(d_{1},d_{2},d_{3}\right)\right)$
into the "matrix>> coordinates
$\left(M\left(m_{1},m_{2},m_{3}\right)\right)$
as follows:
$D\left(d_{1},d_{2},d_{3}\right)\rightarrow M\left(m_{1},m_{2},m_{3}\right)$,
where ${m_{1}=d_{2}}+1$, ${m_{2}=d_{1}}+1$, $m_{3}=d_{3}+1$.

Then it will be necessary to check the Cauchy problem  (\ref{Lyapin:05}) -- (\ref{Lyapin:06}) for solvability, which means, we have to check whether coefficients of the difference operator  (\ref{Lyapin:04}) satisfy to condition  (\ref{Lyapin:06a}).

\section{Example}\label{example}

We consider the polynomial difference operator
\begin{equation}\label{lostformulae}
\begin{array}{c}
P(\delta_{1},\delta_{2},\delta_{3})=c_{000}+c_{100}\delta_{1}+c_{200}\delta_{1}^{2}+c_{010}\delta_{2} + c_{110}\delta_{1}\delta_{2}+c_{210}\delta_{1}^{2}\delta_{2}+\\
+c_{020}\delta_{2}^{2}+c_{120}\delta_{1}\delta_{2}^{2}+c_{220}\delta_{1}^{2}\delta_{2}^{2} + c_{001}\delta_{3}+c_{101}\delta_{1}\delta_{3}+c_{201}\delta_{1}^{2}\delta_{3}+\\
+c_{011}\delta_{2}\delta_{3}+c_{111}\delta_{1}\delta_{2}\delta_{3} + c_{211}\delta_{1}^{2}\delta_{2}\delta_{3}+c_{021}\delta_{2}^{2}\delta_{3}+ c_{121}\delta_{1}\delta_{2}^{2}\delta_{3}+c_{221}\delta_{1}^{2}\delta_{2}^{2}\delta_{3}=\\
=1+2\delta_{1}+3\delta_{1}^{2}+4\delta_{2}+5\delta_{1}\delta_{2}+6\delta_{1}^{2}\delta_{2}+ 1\delta_{2}^{2}+2\delta_{1}\delta_{2}^{2}+\\
+3\delta_{1}^{2}\delta_{2}^{2}+3\delta_{3}+4\delta_{1}\delta_{3}+5\delta_{1}^{2}\delta_{3}+ 6\delta_{2}\delta_{3}+\\
+80\delta_{1}\delta_{2}\delta_{3}+8\delta_{1}^{2}\delta_{2}\delta_{3}+3\delta_{2}^{2}\delta_{3}+ 4\delta_{1}\delta_{2}^{2}\delta_{3}+5\delta_{1}^{2}\delta_{2}^{2}\delta_{3},
\end{array}
\end{equation}
where $b_{x}=2$, $b_{y}=2$, $m=1$. We fix $\beta_{0}=(x_{\beta},y_{\beta},m)=(1,1,1)$, $B_{x}=4$, $B_{y}=3$.
Then the set of initial data is
$$
L=\Pi\setminus\Pi_{\left(1,1\right)},
$$
where
$\Pi=\left\{ (x,y,z)\in\mathbb{Z}^{3},\,0\leqslant x\leqslant4,\,0\leqslant y\leqslant3,z\geqslant0\right\}
$ and $
\Pi_{\left(1,1\right)}=\{(x,y,z)\in\mathbb{Z}_{+}^{3}:1\leqslant x\leqslant3,1\leqslant y\leqslant2,z\geqslant1.\}
$

The three-dimensional matrix $C$ of coefficients of polynomial difference operator \eqref{lostformulae} is given by layers, starting from the lowest one:

\nointerlineskip

\begin{equation}
C =
\left(\begin{array}{c|c}
\begin{array}{c}
\begin{matrix}
c_{000} & c_{100} & c_{200}\\
c_{010} & c_{110} & c_{210}\\
c_{020} & c_{120} & c_{220}
\end{matrix}
\end{array}
&
\begin{array}{c}
\begin{matrix}
c_{001} & c_{101} & c_{201}\\
c_{011} & c_{111} & c_{211}\\
c_{021} & c_{121} & c_{221}
\end{matrix}
\end{array}
\end{array}\right) =
\left(\begin{array}{c|c}
\begin{array}{c}
\begin{matrix}
1 & 2 & 3\\
4 & 5 & 6\\
1 & 2 & 3
\end{matrix}
\end{array}
&
\begin{array}{c}
\begin{matrix}
3 & 4 & 5\\
6 & 80 & 8\\
3 & 4 & 5
\end{matrix}
\end{array}
\end{array}\right).
\end{equation}

\begin{figure}[H]
\includegraphics[width=9 cm]{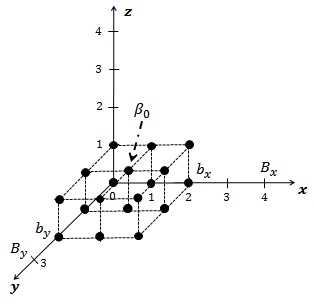}
\caption{The arrangement of the elements of the matrix of coefficients $C$ \label{fig01}}
\end{figure}

The arrangement of the elements of the matrix of coefficients $C$ in the Cartesian coordinate system is shown in Figure~\ref{fig01}.

The problem is to find a value of function
$f\left(x,y,z\right)$ at a point $A$ with coordinates
$\left(2,1,4\right)$.

The matrix $F$ of the initial data is given by layers, starting from the lowest one:

\nointerlineskip

\begin{multline}
F =
\left(\begin{array}{c|c|c|c|c}
F_{0} & F_{1} & F_{2} & F_{3} & F_{4}
\end{array}\right)=
\\
\left(\begin{array}{c|c|c|c|c}
\begin{array}{ccccc}
1 & 2 & 3 & 4 & 5\\
4 & 5 & 6 & 1 & 3\\
1 & 2 & 3 & 8 & 4\\
1 & 2 & 3 & 4 & 6
\end{array}
&
\begin{array}{ccccc}
3 & 4 & 5 & 7 & 6\\
6 & 0 & 0 & 0 & 8\\
3 & 0 & 0 & 0 & 4\\
6 & 5 & 3 & 2 & 3
\end{array}
&
\begin{array}{ccccc}
3 & 4 & 5 & 7 & 6\\
5 & 0 & 0 & 0 & 2\\
3 & 0 & 0 & 0 & 4\\
6 & 8 & 1 & 4 & 3
\end{array}
&
\begin{array}{ccccc}
3 & 4 & 5 & 7 & 6\\
6 & 0 & 0 & 0 & 1\\
3 & 0 & 0 & 0 & 4\\
6 & 6 & 7 & 1 & 4
\end{array}
&
\begin{array}{ccccc}
3 & 4 & 5 & 7 & 6\\
3 & 0 & 0 & 0 & 8\\
3 & 0 & 0 & 0 & 4\\
6 & 2 & 6 & 7 & 3
\end{array}
\end{array}\right).
\end{multline}

\begin{figure}[H]
\includegraphics[width=12.5 cm]{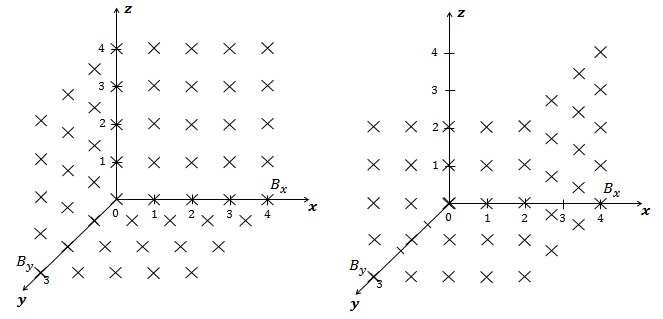}
\caption{The arrangement of the elements of the matrix $F$ \label{fig02}}
\end{figure}

The arrangement of the elements of the matrix $F$ if the initial data in the Cartesian coordinates is shown in Figure~\ref{fig02}.

Now we have to check the following items:
\begin{enumerate}
\item	A transition from the Cartesian coordinates into "matrix>> coordinates:
$\beta_{0}=[1,1,1]\rightarrow\beta_{1}=[2,2,2]$, {$A(2,1,4)\rightarrow A_{1}\left(2,3,5\right)$}.
\item	A solvability of the problem. Since the inequality $$80 = c_{111} > \sum_{\left(\alpha_{1}, \alpha_{2}\right) \neq \left(1,1\right), \alpha_{3} = 1} \|c_{\alpha_{1}\alpha_{2}\alpha_{3}}\| = 38$$
takes place, the problem is solvable.
\item	To find unknown (zero) values in the matrix $F$ on the second layer it is necessary to solve the system of difference equations
\begin{equation}\label{Lyapin:08}
\begin{cases}
P(\delta_{1},\delta_{2},\delta_{3})f(x,y,z)=0,(x,y,z)\in \ensuremath{\Pi}\,for z=0,1,\\
f(x,y,z)=\varphi(x,y,z),(x,y,z)\in L\, for \, z=0,1.
\end{cases}
\end{equation}

The matrix $T$ of the system of difference equations (\ref{Lyapin:08}) will be a block T\"{o}plitz matrix (for details see \cite{ref-book8, ref-journal15}):
\begin{equation}
    T=\left(\begin{array}{cc}
    T_{0} & T_{1}\\
    T_{-1} & T_{0}
    \end{array}\right),
\end{equation}
where
\begin{equation}
    T_{-1}=
    \left(\begin{array}{ccc}
    4 & 5 & 0\\
    3 & 4 & 5\\
    0 & 3 & 4
    \end{array}\right),T_{0}=\left(\begin{array}{ccc}
    80 & 8 & 0\\
    6 & 80 & 8\\
    0 & 6 & 80
    \end{array}\right),T_{1}=\left(\begin{array}{ccc}
    4 & 5 & 0\\
    3 & 4 & 5\\
    0 & 3 & 4
    \end{array}
    \right).
\end{equation}

Solving system (\ref{Lyapin:08}) yields
\begin{equation}
    F_{1}=\left(\begin{array}{ccccc}
    3.0000 & 4.0000 & 5.0000 & 7.0000 & 6.0000\\
    6.0000 & -2.1249 & -1.8591 & -2.6440 & 8.0000\\
    3.0000 & -1.4028 & -1.3917 & -2.4811 & 4.0000\\
    6.0000 & 5.0000 & 3.0000 & 2.0000 & 3.0000
    \end{array}\right).
\end{equation}

\item To find the unknown (zero) values in the matrix $F$ on the third layer, it is necessary to solve the system of difference equations:
\begin{equation}\label{Lyapin:09}
\begin{cases}
P(\delta_{1},\delta_{2},\delta_{3})f(x,y,z)=0,(x,y,z)\in \ensuremath{\Pi} \,for z=1,2,\\
f(x,y,z)=\varphi(x,y,z),(x,y,z)\in L\, for \, z=1,2.
\end{cases}
\end{equation}

The matrix of the system of equations (\ref{Lyapin:09})  will be equal to the matrix of the system of equations (9). Solving the system of equations (\ref{Lyapin:09})  yields:
\begin{equation}
    F_{2}=\left(\begin{array}{ccccc}
    3.0000 & 4.0000 & 5.0000 & 7.0000 & 6.0000\\
    5.0000 & -1.1915 & -0.4974 & -1.9906 & 2.0000\\
    3.0000 & -1.0084 & 0.0246 & -1.2240 & 4.0000\\
    6.0000 & 8.0000 & 1.0000 & 4.0000 & 3.0000
    \end{array}\right).
\end{equation}

\item Continuing the process yields:
\begin{equation}
    F_{3}=\left(\begin{array}{ccccc}
    3.0000 & 4.0000 & 5.0000 & 7.0000 & 6.0000\\
    6.0000 & -1.4035 & -0.7450 & -1.5169 & 1.0000\\
    3.0000 & -1.4785 & -0.3556 & -1.2873 & 4.0000\\
    6.0000 & 6.0000 & 7.0000 & 1.0000 & 4.0000
    \end{array}\right),
\end{equation}

\begin{equation}
    F_{4}=\left(\begin{array}{ccccc}
    3.0000 & 4.0000 & 5.0000 & 7.0000 & 6.0000\\
    3.0000 & -1.1586 & -0.7080 & -2.1118 & 8.0000\\
    3.0000 & -1.1695 & -0.4721 & -1.8794 & 4.0000\\
    6.0000 & 2.0000 & 6.0000 & 7.0000 & 3.0000
    \end{array}\right).
\end{equation}
\end{enumerate}

Thus, the value is $f \left (1,2,4 \right) = - 0.7080$.

\section{Description of the algorithm}

The algorithm was implemented in the MatLab2014 32bit environment. The calculations were performed on an Intel (R) Core (TM) i5-3330S CPU 2.70 GHz, 32bit, 4.00 GB RAM, running Windows 7 Enterprise SP1. The counting time for the given example in Section~\ref{example} was less then a second.

\begin{algorithmic}[1]
\Require a point $\beta_{0}$, a matrix $C$ of coefficients, a point $A$, a matrix $F$ of the initial data
\Ensure the value of the function $f\left(x,y,z\right)$ at a point $A$
\Procedure {DATA}{$C, \beta_{0},A,F$}

\State $\beta_{1}:=$ are coordinates of  $\beta_{0}$ in "matrix" coordinates

\If
    {$\|C\left(\beta_{1}\left(1\right),\beta_{1}\left(2\right), \beta_{1}\left(3\right)\right)\| < \sum\left(\sum\|C\left(:,:, \beta_{1}\left(3\right)\right)\|\right) - \|C\left(\beta_{1}\left(1\right), \beta_{1}\left(2\right),\beta_{1}\left(3\right)\right)\|$}

    \Return {The error of entering the matrix $C$}
\EndIf

\State $A_{1}:=$ are coordinates of $A$ in "matrix" coordinates
\State $Q:=C\left(:,:,end\right)$

\State $p=size\left(Q,2\right)$
\State $NAD=\beta_{1}\left(1\right)-1$
\State $POD=size\left(Q,1\right)-\beta_{1}\left(1\right)$
\State $\begin{array}{c}
sizeTblock=size\left(F,1\right)*size\left(F,2\right)-size\left(F,1\right)*\left(size\left(Q,2\right)-1\right)-\\
-\left(size\left(F,2\right)-\left(size\left(Q,2\right)-1\right)\right)*\left(size\left(Q,1\right)-1\right)
\end{array}$
\State $T:=zeros\left(sizeTblock,sizeTblock\right)$


\State $a=Q\left(\beta_{1}\left(1\right),:\right)$
\State $e=\beta_{1}\left(2\right)$
\State $column_{c}=zeros\left(p,1\right)$
\State $col=1$

\While {$e\geqslant1$}
    \State$column_{c}\left(col\right):=a\left(e\right)$
    \State $col:=col+1$
    \State $e:=e-1$
\EndWhile

\State $e:=\beta_{1}(2)$
\State $row_{c}:=zeros\left(1,p\right)$
\State $r:=1$

\For {$m$ from $e$ to $length\left(a\right)$}
    \State $row_{c}\left(r\right):=a\left(e\right)$
    \State $r:=r+1$
    \State$e:=e+1$
\EndFor

\State $T_{0}:=toeplitz\left(column_{c},row_{c}\right)$


\For {$i$ from $0$ to $sizeTblock/(size(Q,2))-1$}
    \State $T(((i*p)+1):p*(i+1),((i*p)+1):p*(i+1)):=T_{0}$
\EndFor


\For {$w$ from $1$ to $NAD$}
    \State $a:=Q\left(\beta_{1}\left(1\right)-w,:\right)$
    \State $e:=\beta_{1}\left(2\right)$
    \State $column_{c}:=zeros\left(p,1\right)$
    \State $col:=1$

    \While {$e\geqslant1$}
        \State $column_{c}\left(col\right):=a\left(e\right)$
        \State $col:=col+1$
        \State $e:=e-1$
    \EndWhile

    \State $e:=\beta_{1}\left(2\right)$
    \State $row_{c}=zeros\left(1,p\right)$
    \State $r:=1$

    \For {$m$ \textbf{from} $e$ \textbf{to} $length\left(a\right)$}
        \State $row_{c}\left(r\right):=a\left(e\right)$
        \State $r:=r+1$
        \State $e:=e+1$
    \EndFor

    \State $T_{NAD}:=toeplitz\left(colomn_{c},row_{c}\right)$


    \For {$i$ from $0$ to $sizeTblock/\left(size\left(Q,2\right)\right)-1-w$}
        \State $T((((i*p)+1)):p*(i+1),((i*p)+1)+p*w:p*(i+1)+p*w):=T_{NAD}$
    \EndFor
\EndFor


\For {$w$ from $1$ to $POD$}
    \State $a:=Q\left(\beta_{1}\left(1\right)+w,:\right)$
    \State $e:=\beta_{1}\left(2\right)$
    \State $column_{c}=zeros\left(p,1\right)$;
    \State $col=1$;

    \While {$e\leqslant1$}
        \State $column_{c}\left(col\right)=a\left(e\right)$
        \State $col=col+1$
        \State $e=e-1$
    \EndWhile

    \State $e=\beta_{1}\left(2\right)$
    \State $row_{c}=zeros\left(1,p\right)$
    \State $r=1$

    \For {$m$ from $e$ to $length\left(a\right)$}
        \State $row_{c}\left(r\right)=a\left(e\right)$
        \State $r=r+1$
        \State $e=e+1$
    \EndFor

    \State $T_{POD}=toeplitz\left(column_{c},row_{c}\right)$


    \For {$i$ from $0$ to $sizeTblock/\left(size\left(Q,2\right)\right)-1-w$}
        \State $T(((i*p)+1)+p*w:p*(i+1)+p*w,((i*p)+1):p*(i+1)):=T_{POD}$
    \EndFor
\EndFor


\State $b:=zeros\left(size\left(T,1\right),1\right)$;

\For {$k$ from $\beta_{1}\left(3\right)$ to $A\left(3\right)$}
    \State $t:=1$;
    \For {$i$ from $0$ to $\left(size\left(F,2\right)-size\left(C,2\right)\right)$}
        \For { $j$ from $0$ to $\left(size\left(F,1\right)-size\left(C,1\right)\right)$}
            \State $f_{iter} := F\left(j + 1:size\left(C,1\right) + j,\left(i + 1\right):size\left(C,2\right) + i,k - size\left(C,3\right) + 1:k\right)$
            \State $f_{new}:=f_{iter}.*C$
            \State $b\left(t\right)=-\left(\sum\left(\sum\left(\sum\left(f_{new}\right)\right)\right)\right)$
            \State $t:=t+1$
        \EndFor
        \State $elements:=linsolve\left(T,b\right)$
        \State $m:=1$

        \For {$i$ from $0$ to $\left(size\left(F,2\right)-size\left(C,2\right)\right)$}
            \For {$j$ from $0$ to $\left(size\left(F,1\right)-size\left(C,1\right)\right)$}
                \State $F\left(\beta_{1}\left(1\right) + j, \beta_{1}\left(2\right)+i, k\right) := elements\left(m\right)$;
                \State $m:=m+1$
            \EndFor
        \EndFor
    \EndFor
\EndFor

\Return {$f\left(A\right)$}

\EndProcedure
\end{algorithmic}

The complete code of the program is available at http://github.com/ApanovichMS/CauchyParallelepiped.git.

\section{Conclusion}
Classical methods for solving differential equations (methods of Runge-Kutta of the fourth order, Euler, Newton, etc.) have certain difficulties in their application in the multidimensional case. In each specific task, it is necessary to individually select a solving method. This makes the development of a universal approach to solving multidimensional differential equations using these methods laborious. In this study, we expanded the possibility of using standard symmetric difference schemes to approximate two-dimensional differential equations with constant coefficients, where the problem is considered "on the plane>>, in which, as a rule, time and one spatial variable are used as independent variables. The scheme proposed by us, firstly, has an arbitrary number of points, and secondly, it can be used to approximate differential equations in the three-dimensional case. For the heat conduction equation, this would mean that when describing heat transfer, it becomes possible to take into account the anisotropy of this process by adding an additional spatial variable.

The process of approximating differential equations by a difference scheme in the three-dimensional case gives a system of difference equations, the matrix of which has a block T\"{o}plitz form. The process of solving such systems for a large number of variables and initial data is associated with high computational costs. In our work, we proposed an algorithm for solving such a system using values of the
coefficients of the difference equation, which made it possible to automate the solution of the Cauchy problem with the initial data in the parallelepiped. Thus, our proposed approach defines a unified algorithm for solving differential equations when they are approximated by a difference scheme in a parallelepiped.


\section*{Funding}
This work is supported by the Krasnoyarsk Mathematical Center and financed by the Ministry of
Science and Higher Education of the Russian Federation in the framework of the establishment
and development of regional Centers for Mathematics Research and Education (Agreement No.
075-02-2021-1388).

\end{document}